\documentclass[12pt]{article}
\usepackage{amsmath}
\usepackage{amssymb}
\newcommand{\R}{{\mathbb R}}
\newcommand{\N}{{\mathbb N}}
\newcommand{\fn}{\!:\!}
\newcommand{\hk}{{\cal HK}}
\newcommand{\lsum}{\sum\limits}
\newcommand{\lsup}{\sup\limits}
\newcommand{\qed}{\mbox{$\quad\blacksquare$}}
\newtheorem{theorem}{Theorem}
\newtheorem{lemma}[theorem]{Lemma}

\newtheorem{example}[theorem]{Example}
\begin{document}
\hspace{-2cm}
\raisebox{12ex}[1ex]{\fbox{{\footnotesize
To appear in {\it Czechoslovak Mathematical Journal.}
Accepted November~3, 2003
}}}
$~$ \\ [.5in]
\begin{center}
{\large\bf Estimates of the remainder in Taylor's theorem using the
Henstock--Kurzweil integral}
\vskip.25in
Erik Talvila\footnote{Research partially supported by the
Natural Sciences and Engineering Research Council of Canada.
An adjunct appointment in the Department of Mathematical and
Statistical Sciences, University of Alberta, made valuable library
and computer resources available.}\\ [2mm]
{\footnotesize
Department of Mathematics and Statistics \\
University College of the Fraser Valley\\
Abbotsford, BC Canada V2S 7M8\\
Erik.Talvila@ucfv.ca}
\end{center}

{\footnotesize
\noindent
{\bf Abstract.} When a real-valued function of one variable is
approximated by its $n^{th}$ degree Taylor polynomial, the
remainder is estimated using the Alexiewicz and Lebesgue $p$-norms in cases
where $f^{(n)}$ or $f^{(n+1)}$ are Henstock--Kurzweil integrable.
When the only assumption is that $f^{(n)}$ is Henstock--Kurzweil integrable
then a modified form of the $n^{th}$ degree Taylor polynomial is used.
When the only assumption is that $f^{(n)}\in C^0$ then the remainder
is estimated by applying the Alexiewicz norm to Schwartz
distributions of order 1.
}
\vskip.25in

\section{Introduction}
In this paper the Henstock--Kurzweil integral is used to give various
estimates of the remainder in Taylor's theorem in terms of Alexiewicz
and Lebesgue $p$-norms.  Let $[a,b]$ be a compact interval in $\R$ and let
$f\fn[a,b]\to\R$.  Let $n$ be a positive integer.
When $f$ is
approximated by its $n^{th}$ degree Taylor polynomial about $a$, the
remainder is estimated using the Alexiewicz norm of $f^{(n+1)}$ and 
$p$-norms of $f^{(n)}$ in the case
when $f^{(n+1)}$ is Henstock--Kurzweil integrable (Theorem~\ref{n}).  
When the only assumption is that $f^{(n)}$ is Henstock--Kurzweil integrable
then $f^{(n)}$ need not exist at $a$.  In this case we use a modified 
form of the Taylor polynomial where $f^{(n)}$ is evaluated at a point
$x_0\in[a,b]$ at which $f^{(n)}$ exists.  
The resulting modified remainder is then estimated in
the Alexiewicz and $p$-norms (Theorem~\ref{n-1}).
And, when the only assumption is that $f^{(n)}\in C^0$, the 
Alexiewicz norm in used in the space of Schwartz distributions of order 1
to estimate the remainder (Theorem~\ref{C}).
The results extend those in \cite{anastassiou}.
See also \cite{vyborny} for another form of the remainder.

For the Henstock--Kurzweil integral
we have the following version of Taylor's theorem.
\begin{theorem}\label{taylor}
Let $f\fn[a,b]\to\R$ and let $n$ be a positive integer.  If $f^{(n)}\in
ACG_*$ then for all $x\in[a,b]$ we have $f(x) =P_n(x) + R_n(x)$ where
\begin{equation}
P_n(x)=\lsum_{k=0}^n\frac{f^{(k)}(a)(x-a)^k}{k!}\label{polynomial}
\end{equation}
and
\begin{equation}
R_n(x)=\frac{1}{n!}\int_a^x f^{(n+1)}(t)(x-t)^{n}\,dt.\label{remainder}
\end{equation}
\end{theorem}

A proof under
the assumption that $f^{(n)}$ is continuous on $[a,b]$
and $f^{(n+1)}$ exists nearly
everywhere on $(a,b)$ is given in \cite{thompson}.
The general case follows from the Fundamental Theorem:
If $F\fn[a,b]\to\R$ and $F\in ACG_*$ then $F'$ exists almost everywhere,
$F'$ is Henstock--Kurzweil integrable and $\int_a^x F'=F(x) - F(a)$ for
all $x\in [a,b]$.  For the wide Denjoy integral the corresponding function space is $ACG$.
If $F\in ACG$ then $\int_a^x F'_{{\rm ap}}=F(x)-F(a)$ for
all $x\in [a,b]$.  Here the integral is the wide Denjoy integral and the
approximate derivative is used.  All the results proved in the paper  
have a suitable extension to the wide Denjoy integral.
The function spaces $ACG_*$ and $ACG$ are defined in
\cite{celidze} and \cite{saks}.  Note that $AC\subsetneq ACG_*\subsetneq ACG
\subsetneq C^0$.
The set of continuous  functions  that are differentiable nearly everywhere
is a proper subset of $ACG_*$, which is itself a proper subset of the set 
of continuous functions that are differentiable almost everywhere.
The other half of the Fundamental Theorem says that if $f\in\hk$ then
$\frac{d}{dx}\int_a^xf=f(x)$ almost everywhere (and certainly at points
of continuity of $f$).

The Alexiewicz norm of an integrable function $f$ is defined
\begin{equation}
\|f\|=\lsup_{a\leq x\leq b} \left|\int_a^x f\right|.
\end{equation}
We will write $\|f\|_I$ when it needs to be made clear the norm is
over interval $I$.
See \cite{swartz} for a discussion of the Alexiewicz norm and the 
Henstock--Kurzweil
integral.  An equivalent norm is 
$\sup_{(c,d)\subset(a,b)} \left|\int_c^d f\right|$.
It leads to similar estimates to those obtained for the 
Alexiewicz norm in 
Theorems~\ref{n-1}, 
\ref{n} and \ref{C}.  

Denote the space of Henstock--Kurzweil integrable functions on $[a,b]$
by $\hk$.     
Note that if $f^{(n)}\in ACG_*$ then $f^{(n+1)}\in\hk$.

\section{Estimates when $f^{(n-1)}\in ACG_*$}
If $f^{(n-1)}\in ACG_*$ (i.e., $f^{(n)}\in\hk$)
then $f^{(n)}$ need only  exist almost everywhere on
$(a,b)$.  If $f^{(n)}$ exists at $x_0\in[a,b]$ then we  can
modify the Taylor polynomial 
\eqref{polynomial} so that the $n^{th}$ derivative is evaluated at $x_0$
and then
obtain estimates on the resulting remainder term.  The following
lemma is used.
\begin{lemma}\label{lemma1}
If $f\fn[a,b]\to\R$ and $f^{(n-1)}\in ACG_*$ then let $x_0\in[a,b]$ such
that $f^{(n)}$ exists at $x_0$.  Define the modified Taylor polynomial by
\begin{equation}
P_{n,x_0}(x)=P_{n-1}(x)+\frac{f^{(n)}(x_0)(x-a)^n}{n!}
\end{equation}
and define the modified remainder by $R_{n,x_0}(x)=f(x)-P_{n,x_0}(x)$.
Then for all $x\in[a,b]$
\begin{equation}
R_{n,x_0}(x)=
\frac{1}{(n-1)!}\int_a^x\left[f^{(n)}(t)-f^{(n)}(x_0)\right](x-t)^{n-1}
\,dt.\label{remainder2}
\end{equation}
\end{lemma}
The proof is the same as for Lemma~1 in \cite{anastassiou} (which is false
without the proviso that $f^{(n)}$ exists at $x_0$).  Of course
we can take $x_0$ as close to $a$ as we like.  See \cite{folland} for
remainder estimates based on estimates of $|f^{(n)}(t)-f^{(n)}(x_0)|$.

Theorem~2 in \cite{anastassiou} gives pointwise estimates 
of $R_n$ in terms of $p$-norms of $f^{(n)}$ when $f^{(n-1)}\in AC$ and
$f^{(n)} \in L^p$ ($1\leq p\leq \infty$).  We have the following
analogue when $f^{(n-1)}\in ACG_*$.  Note that if $f^{(n-1)}\in ACG_*\setminus
AC$ (i.e., $f^{(n)}\in\hk\setminus L^1$)
then for each $1\leq p\leq\infty$ we have $f^{(n)} \not\in L^p$.
However, we can use the Alexiewicz norm to estimate $R_{n,x_0}$.
\begin{theorem}\label{n-1}
With the notation and assumptions of Lemma~\ref{lemma1},
\begin{equation}
\|R_{n,x_0}\|\leq \frac{(b-a)^n}{n!}\|f^{(n)}(\cdot)
-f^{(n)}(x_0)\|.
\label{n-1.1}
\end{equation}
For all
$x\in[a,b]$,
\begin{equation}
|R_{n,x_0}(x)|\leq \frac{(x-a)^{n-1}}{(n-1)!}\|f^{(n)}(\cdot)-f^{(n)}(x_0)\|
_{[a,x]}.
\label{n-1.2}
\end{equation}
And,
\begin{equation}
\|R_{n,x_0}\|_p\leq\left\{\begin{array}{cl}
\frac{(b-a)^{n-1+1/p}}{(n-1)!\,[(n-1)p+1]^{1/p}}\|f^{(n)}(\cdot)-f^{(n)}(x_0)\|,
& 1\leq p<\infty\\
\frac{(b-a)^{n-1}}{(n-1)!}\|f^{(n)}(\cdot)-f^{(n)}(x_0)\|,
& p=\infty.
\end{array}
\right.
\label{n-1.3}
\end{equation}
\end{theorem}

\noindent
{\bf Proof:} Let $a\leq c\leq b$.  Using Lemma~\ref{lemma1} and the reversal
of integrals criterion in \cite[Theorem~57, p.~58]{celidze},
\begin{eqnarray}
\int_a^c R_{n,x_0} & = & \frac{1}{(n-1)!}\int_a^c\left[
f^{(n)}(t)-f^{(n)}(x_0)\right]\int_t^c(x-t)^{n-1}\,dx\,dt\notag\\
 & = & \frac{(c-a)^n}{n!}\!\!\int_a^\xi\left[\!
f^{(n)}(t)-f^{(n)}(x_0)\!\right]dt\,\,\,\,\,\text{ for some } 
\xi\in [a,c].\label{mvt1}
\end{eqnarray}
Equation \eqref{mvt1} comes from the second mean value theorem for
integrals \cite{celidze}.  Taking the supremum over $c\in[a,b]$ now
gives \eqref{n-1.1}.

Similarly,
\begin{eqnarray}
|R_{n,x_0}(x)| & = & \frac{(x-a)^{n-1}}{(n-1)!} \left|\int_a^\xi\left[
f^{(n)}(t)-f^{(n)}(x_0)\right]\,dt\right|\quad\text{for some } \xi\in[a,x]
\notag\\
 & \leq & \frac{(x-a)^{n-1}}{(n-1)!}\|f^{(n)}(\cdot)-f^{(n)}(x_0)\|.
\end{eqnarray}
This gives \eqref{n-1.2}.  The other estimates in \eqref{n-1.3} 
follow from this.\qed

An alternative approach in Theorem~\ref{n-1} is to assume $f^{(n-1)}$
is $ACG_*$ on $[a,b]$ and $f^{(n)}$ exists at $a$.  Then
equation \eqref{n-1.1}
is replaced with $\|R_n\|\leq (b-a)^n\|f^{(n)}(\cdot)-f^{(n)}(a)\|/n!$
with similar changes in \eqref{n-1.2} and \eqref{n-1.3}.
This is Young's expansion theorem \cite[p.~16]{young}
but with a more precise error estimate.
\section{Estimates when $f^{(n)}\in ACG_*$}
Corollary~1 in \cite{anastassiou} gives a pointwise estimate of $R_n$
in terms of $p$-norms of $f^{(n+1)}$ when $f^{(n)}\in AC$ and 
$f^{(n+1)}\in L^p$.  When $f^{(n)}\in ACG_*\setminus AC$ 
(i.e., $f^{(n+1)}\in\hk\setminus L^1$) then for no
$1\leq p\leq\infty$ do we have $f^{(n+1)}\in L^p$.  But, we can estimate
$R_n$ using the Alexiewicz norm of $f^{(n+1)}$ and $p$-norms of
$f^{(n)}$.
\begin{theorem}\label{n}
If $f\fn[a,b]\to\R$ such that $f^{(n)}\in ACG_*$ then
\begin{equation}
\|R_n\|\leq \frac{(b-a)^{n+1}}{(n+1)!}\|f^{(n+1)}\|.
\label{n.1}
\end{equation}
For all
$x\in[a,b]$,
\begin{equation}
|R_n(x)|\leq \frac{(x-a)^{n}}{n!}\|f^{(n+1)}\|_{[a,x]}.
\label{n.2}
\end{equation}
And,
\begin{equation}
\|R_n\|_p\leq\left\{\begin{array}{cl}
\frac{(b-a)^{n+1/p}}{n!\,(np+1)^{1/p}}\|f^{(n+1)}\|,
& 1\leq p<\infty\\
\frac{(b-a)^{n}}{n!}\|f^{(n+1)}\|,
& p=\infty.
\end{array}
\right.
\label{n.3}
\end{equation}
Also,
\begin{equation}
\|R_n\|_p\leq\left\{\begin{array}{cl}
\frac{(b-a)^{n+1/p}}{n!\,(np+1)^{1/p}}|f^{(n)}(a)| +A_i,
& 1\leq p<\infty\\
\frac{(b-a)^{n}}{n!}|f^{(n)}(a)| +\frac{(b-a)^{n}}{n!}\|f^{(n)}\|_\infty,
& p=\infty,
\end{array}
\right.
\label{n.4}
\end{equation}
where $1/\alpha+1/\beta=1$ and
\begin{eqnarray*}
A_1 & = & \frac{(b-a)^{n-1+1/p}\|f^{(n)}\|}{(n-1)!\,[(n-1)p+1]^{1/p}}, \quad
\text{for } n\geq 1\\
A_1 & = & \|f(\cdot)-f(a)\|_p, \quad
\text{for } n= 1\\
A_2 & = & \frac{(b-a)^{n-1+1/p+1/\beta}\|f^{(n)}\|_\alpha}{(n-1)!\,
[(n-1)\beta +1]^{1/\beta}[(n-1+1/\beta)p+1]^{1/p}}, \quad
\text{for } 1<\alpha\leq\infty\\
A_2 & = & \frac{(b-a)^{n-1+1/p}\|f^{(n)}\|_1}{(n-1)!\,
[(n-1)p+1]^{1/p}}, \quad
\text{for } \alpha=1\\
A_3 & = & \frac{(b-a)^{1-1/p}}{(n-1)!\,
[(n-1)p+1]^{1/p}}
\left(\int_a^b|f^{(n)}(t)|^p(b-t)^{(n-1)p+1}\,dt\right)^{1/p}\\
A_4 & = & \frac{(b-a)^{n(1-1/p)}}{n!}
\left(\int_a^b|f^{(n)}(t)|^p(b-t)^{n}\,dt\right)^{1/p}.
\end{eqnarray*}
\end{theorem}

\noindent
{\bf Proof:}
The proof of \eqref{n.1} is very similar to the proof of \eqref{n-1.1},
except that we begin with the remainder in the form of \eqref{remainder}.

Using the second mean value theorem, we have
\begin{equation}
|R_n(x)|  =  \frac{(x-a)^n}{n!}\left|\int_a^\xi f^{(n+1)}(t)\,dt\right|
\quad \text{for some } \xi\in[a,x].
\end{equation}
The estimates in \eqref{n.2} and \eqref{n.3} now follow directly.

Integrate \eqref{remainder} by parts to get
\begin{equation}
R_n(x) = -\frac{f^{(n)}(a)(x-a)^n}{n!} + \frac{1}{(n-1)!}
\int_a^x f^{(n)}(t)(x-t)^{n-1}\,dt.\label{altremainder}
\end{equation}
Then
\begin{equation}
\|R_n\|_p \leq \frac{|f^{(n)}(a)|}{n!}I_1^{1/p} + \frac{1}{(n-1)!}I_2^{1/p}
\end{equation}
where
\begin{equation}
I_1=\int_a^b(x-a)^{np}\,dx=\frac{(b-a)^{np+1}}{np+1}
\end{equation}
and
\begin{equation}
I_2=\int_a^b\left|\int_a^xf^{(n)}(t)(x-t)^{n-1}\,dt\right|^p\,dx.
\end{equation}
$A_1$ is obtained from $I_2$ using the second mean value theorem
and $A_2$ using H\"{o}lder's inequality.

Writing 
\begin{equation}
I_2\leq \int_a^b\left(\int_a^x|f^{(n)}(t)|(x-t)^{n-1}\frac{dt}{
x-a}\right)^p(x-a)^p\,dx,
\end{equation}
Jensen's inequality and Fubini's theorem give
\begin{eqnarray}
I_2 & \leq & \int_a^b|f^{(n)}(t)|^p\int_t^b(x-t)^{(n-1)p}(x-a)^{p-1}dx\,dt
\label{altnorm}\\
 & \leq & \frac{(b-a)^{p-1}}{(n-1)p+1}\int_a^b|f^{(n)}(t)|^p(b-t)^{(n-1)p+1}
dt.
\end{eqnarray}
From this we obtain $A_3$.

Using Jensen's inequality in the form
\begin{eqnarray}
I_2 & \leq & \int_a^b\left(\int_a^x|f^{(n)}(t)|\frac{(x-t)^{n-1}n\,dt}{
(x-a)^n}\right)^p\frac{(x-a)^{np}}{n^p}dx\\
 & \leq & \frac{1}{n^{p-1}}\int_a^b\int_a^x |f^{(n)}(t)|^p(x-t)^{n-1}
(x-a)^{n(p-1)}dt\,dx,
\end{eqnarray}
we can apply Fubini's theorem to get
\begin{equation}
I_2\leq \frac{(b-a)^{n(p-1)}}{n^p}\int_a^b|f^{(n)}(t)|^p
(b-t)^n\,dt,
\end{equation}
which gives $A_4$.  The case $p=\infty$ follows directly 
from \eqref{altremainder}.\qed

Note that $A_2$, $A_3$ and $A_4$ all lead to estimates of form
$A_k\leq C_{n,p}(b-a)^n\|f^{(n)}\|_p$ ($k=2,3,4$) where $C_{n,p}$
is independent of $a$, $b$ and $f$.

The integral over $x$ in \eqref{altnorm} can be evaluated using
hypergeometric functions.  However, this does not markedly improve
the estimate for $A_3$.  Similarly with $A_4$.

In \eqref{n.2} we have
$R_n(x) =o[(x-a)^n]$ as $x\to a$.

Note that if $f^{(n)}\in ACG_*$ then $\|f^{(n+1)}\|=\max_{a\leq x\leq b}
|f^{(n)}(x)-f^{(n)}(a)|$.  This affects  how the remainder is written
in Theorems~\ref{n-1} and \ref{n}.

\begin{example}
{\rm
Let $0<a_n\leq 1$ be a sequence that decreases to $0$.  Let $b_n$
be a positive sequence that decreases to 0 such that the intervals
$(a_n-b_n,a_n+b_n)$ are disjoint.  For this it suffices to take
$b_n\leq\min([a_{n-1}-a_n]/2,[a_n-a_{n+1}]/2)$.  Let
$f_n(x)=(x-a_n+b_n)^2(x-a_n-b_n)^2$ for $|x-a_n|\leq b_n$ and 0, otherwise.
Let $\alpha>0$ and define $f(x)=\sum n^{\alpha}f_n(x)$.
For $|x-a_n|< b_n$ we have $f_n'(x)=4(x-a_n+b_n)(x-a_n-b_n)(x-a_n)$
and $f''_n(x)=4[3(x-a_n)^2-b_n^2]$.
Suppose
$b_n=o(a_n)$ as $n\to\infty$ and $n^\alpha b_n^3\to 0$.  Then
$\max_{|x-a_n|\leq b_n}|f_n'(x)|=O(b_n^3)$ so
$f\in C^1_{[0,1]}$ but $f\not\in C^2_{[0,1]}$.  If $\sum n^\alpha b_n^3
=\infty$ then $f''\in\hk\setminus L^1$. Let $a=0$ in Theorems~\ref{n-1}
and \ref{n}.

(i) If $n^\alpha b_n^3/a_n\not\to 0$ then $f''(0)$ does not exist.  
An example of this
case is $a_n=1/n$, $b_n=c\,n^{-\beta}$ for $\beta\geq 2$ 
and small enough $c$, and 
$3\beta -1\leq\alpha<3\beta$.  Let 
$x_0\in(0,1)\setminus\{a_n\pm b_n\}_{n\in\N}$.  Then
$f''(x_0)$ exists and we have the modified second degree Taylor polynomial
$P_{2,x_0}(x)=f(0)+f'(0)x+f''(x_0)x/2=f''(x_0)x/2$.
The modified remainder is $R_{2,x_0}(x)=\int_0^x[f''(t)-f''(x_0)](x-t)\,dt$.
If we take $x_0=a_n\pm b_n/\sqrt{3}$ for some $n$ then $f''(x_0)=0$ and
we have $\|f''(\cdot)-f''(x_0)\|_{[0,x]}\leq\max_{0\leq 
y\leq 1/\lfloor{1/x}\rfloor}|f'(y)|
=\sup_{n\geq \lfloor{1/x}\rfloor}n^{\alpha}|f_n'(a_n\pm b_n/\sqrt{3})|=
(8c^3/3\sqrt{3})\sup_{n\geq \lfloor{1/x}\rfloor}n^{\alpha-3\beta}=
(8c^3/3\sqrt{3}){\lfloor{1/x}\rfloor}^{\alpha -3\beta}$.  
This allows us to obtain all of
the estimates in Theorem~\ref{n-1}.  Note that $|R_{2,x_0}(x)|
=O(x{\lfloor{1/x}\rfloor}^{\alpha -3\beta})=O(x^{1-\alpha +3\beta})$
as $x\to 0$ (and $1<1-\alpha +3\beta\leq 2$).

(ii) If $n^\alpha b_n^3/a_n\to 0$ then $f''(0)$ exists.
An example is $0< \alpha<3\beta-1$ with $a_n$ and  $b_n$ as in (i).
This gives $f''\in L^1$.
Now $P_{2,0}(x)=0$ and $R_{2,0}(x)=\int_0^xf''(t)(x-t)\,dt$.  The
quantity $\|f''\|_{[0,x]}$ can be estimated as in (i).

(iii) We can use Theorem~\ref{n} with $n=1$ to get $P_1(x)=0$ 
and $R_1(x)=\int_0^xf''(t)(x-t)\,dt$.  This leads to the same estimates
for $\|f''\|_{[0,x]}$ as in (i) in the case $x_0=a_n\pm b_n/\sqrt{3}$.
\qed
}
\end{example}
\section{Estimates when $f^{(n)}\in C^0$}
We now show that \eqref{n.1} continues to hold when the only assumption
is that $f^{(n)}\in C^0$.  Under the Alexiewicz norm, the space of
Henstock--Kurzweil integrable functions  is not 
complete.  It's completion is the subspace of distributions that
are the distributional derivative of a continuous function, i.e.,
distributions of order 1 (see
\cite{swartz}).  
Thus, if $f$ is in the completion of $\hk$ then $f\in{\cal D}'$ (Schwartz
distributions) and
there is a continuous
function $F$, vanishing at $a$, such that $F'(\phi)=-F(\phi')=
-\int_{a}^b F\phi'=f(\phi)$
for all test functions $\phi\in{\cal D}=\{\phi\fn[a,b]\to\R\mid
\phi\in C^\infty \text{ and } {\rm supp}(\phi)\subset(a,b)\}$.  
And, we can
compute the Alexiewicz norm of $f$ via $\|f\|=\max_{x\in[a,b]}|F(x)|=
\|F\|_\infty$.
\begin{theorem}\label{C}
If $f\fn[a,b]\to\R$ such that $f^{(n)}\in C^0$ then
\begin{equation}
\|R_n\|\leq \frac{(b-a)^{n+1}}{(n+1)!}\|f^{(n+1)}\|=\frac{(b-a)^{n+1}}{(n+1)!}
\max_{a\leq x\leq b}|f^{(n)}(x)-f^{(n)}(a)|.
\label{C.1}
\end{equation}
\end{theorem}

\noindent
{\bf Proof:}
From Lemma~\ref{lemma1}, 
$
R_{n}(x)=
\frac{1}{(n-1)!}\int_a^x\left[f^{(n)}(t)-f^{(n)}(a)\right](x-t)^{n-1}
\,dt$.  Integrating $x$ from $a$  to $y$ and reversing orders of 
integration gives
$\int_a^yR_n=\frac{1}{n!}\int_a^y\left[f^{(n)}(t)-f^{(n)}(a)\right](y-t)^{n}
\,dt$.  Equation \eqref{C.1} now follows.
\qed
\begin{example}
{\rm
Let $g$ be continuous and nowhere differentiable on $[0,1]$.  Let $n\geq 1$
and define $f(x)=\frac{1}{(n-1)!}\int_0^xg(t)(x-t)^{n-1}\,dt$.  By
differentiating under the integral and then using the Fundamental Theorem,
$f^{(n)}=g\in C^0$ but $f^{(n)}\not\in ACG_*$.  We have 
$R_n(x)=\frac{1}{(n-1)!}\int_0^x[g(t)-g(0)](x-t)^{n-1}\,dt$ and
$\|R_n\|\leq \frac{1}{(n+1)!}\max_{0\leq x\leq 1}|g(x)-g(0)|$.\qed
}
\end{example}

\end{document}